\pdfoutput=1 
\documentclass[conference, 10pt]{IEEEtran}

\usepackage{todonotes}

\usepackage{ifthen}
\usepackage{graphicx}
\graphicspath{{./img/}}

\usepackage{booktabs} 
\usepackage{cite} 
\usepackage[cmex10]{amsmath}
\usepackage{amssymb} 
\usepackage{amsfonts}
\usepackage{amsthm}

\usepackage{fixltx2e}

\usepackage{algorithm} \usepackage{algorithmic}

\usepackage{subfigure}

\newcommand{\eg}{{\it e.g.}}  \newcommand{\ie}{{\it i.e.}}

 \newcommand{\reals}{{\mathbb R}}
 
\newcommand\plusvee{\mathrel{\ooalign{\lower.5ex
\hbox{$\scriptstyle\vee\mkern.5mu$}\cr\hidewidth\raise.450ex
\hbox{$\scriptstyle+$}\cr}}}

\newcommand{\argmin}{\operatornamewithlimits{argmin}}

\newcommand{\minimize}{\operatornamewithlimits{minimize}}

\newboolean{draft}

\newcommand{\isdraft}[2]{\ifthenelse{\boolean{draft}}{#1}{#2}}

\newtheorem{theorem}{Theorem}

\newtheorem{remark}[theorem]{Remark}

\title{Distributed, simple and stable network localization}

\author{\IEEEauthorblockN{Cl\'audia Soares*, Jo\~ao Xavier, and Jo\~ao
    Gomes}
  \IEEEauthorblockA{
    Institute for Systems and Robotics (ISR), Instituto Superior
    T\'ecnico, Universidade de
    Lisboa\\
    Lisbon, Portugal\\
    \{csoares,jxavier,jpg\}@isr.ist.utl.pt\\}
\thanks{ This research was
    supported by Fundação para a Ci\^encia e a Tecnologia (projects
    PEst-OE/EEI/LA0009/2013, CMU-PT/SIA/0026/2009 and PhD grant
    SFRH/BD/72521/2010) and EU FP7 project MORPH (grant agreement no.\
    288704)}}

\begin{document} 

\maketitle
\thispagestyle{empty}
\pagestyle{empty}

\begin{abstract}
  We propose a simple, stable and distributed algorithm which directly
  optimizes the nonconvex maximum likelihood criterion for sensor
  network localization, with no need to tune any free parameter. We
  reformulate the problem to obtain a gradient Lipschitz cost; by
  shifting to this cost function we enable a Majorization-Minimization
  (MM) approach based on quadratic upper bounds that decouple across
  nodes; the resulting algorithm happens to be distributed, with all
  nodes working in parallel. Our method inherits the MM stability:
  each communication cuts down the cost function. Numerical
  simulations indicate that the proposed approach tops the performance
  of the state of the art algorithm, both in accuracy and
  communication cost.
\end{abstract}



\section{Problem statement}
\label{sec:problem-statement} 
The sensor network is represented as an undirected connected graph
$\mathcal{G} = (\mathcal{V},\mathcal{E})$. The node set $\mathcal{V} =
\{1,2, \dots, n\}$ denotes the sensors with unknown positions. There
is an edge $i \sim j \in {\mathcal E}$ between sensors $i$ and $j$ if
a noisy range measurement between nodes $i$ and $j$ is available at
both, and if $i$ and $j$ can communicate with each other.  The
set of sensors with known positions, named anchors, is denoted by
${\mathcal A} = \{ 1, \ldots, m \}$. For each $i \in {\mathcal V}$, we
let ${\mathcal A}_i \subset {\mathcal A}$ be the subset of anchors (if
any) relative to which node $i$ also possesses a noisy range
measurement.

Let $\reals^p$ be the space of interest ($p=2$ for planar networks,
and $p=3$ otherwise), $x_i \in \reals^p$ the position of sensor $i$,
and $d_{ij}$ the noisy range measurement between sensors $i$ and $j$,
known by both $i$ and $j$. Without loss of generality, we assume
$d_{ij} = d_{ji}$. Anchor positions are denoted by $a_{k} \in
\reals^{p}$. Similarly, $r_{ik}$ is the noisy range measurement
between sensor $i$ and anchor $k$, available at sensor $i$.

The distributed network localization problem addressed in this work
consists in estimating the sensors' positions $x = \{ x_i\, : \, i \in
\mathcal{V} \}$, from the available measurements $\{ d_{ij} \, : \, i
\sim j \} \cup \{ r_{ik} \, : \, i \in {\mathcal V}, k \in {\mathcal
A}_i \}$, through collaborative message passing between neighboring
sensors in the communication graph~${\mathcal G}$.

Under the assumption of zero-mean, independent and
identically-distributed, additive Gaussian measurement noise, the
maximum likelihood estimator for the sensor positions is the solution
of the optimization problem
\begin{equation}
  \label{eq:snlOptimizationProblem} 
  \minimize_{x} f(x),
\end{equation} 
where
\begin{equation*} 
  f(x) = \sum _{i \sim j} \frac{1}{2}(\|x_{i} - x_{j}\| - d_{ij})^2 + \sum_{i}
  \sum_{k \in \mathcal{A}_{i}} \frac{1}{2}(\|x_{i}-a_{k}\| - r_{ik})^2.
\end{equation*} 
Problem~\eqref{eq:snlOptimizationProblem} is nonconvex and difficult
to solve.  Even in the centralized setting (\ie, all measurements are
available at a central node) currently available iterative techniques
don't claim convergence to the global optimum. Also, even with
noiseless measurements, multiple solutions might exist due to
ambiguities in the network topology itself\cite{AndersonShamesMaoFidan2010}.

\paragraph*{Related work}
\label{sec:related-work}
The literature on sensor network localization can be divided on
centralized and distributed approaches; the first category encompasses
methods that use a central processing node, which collects all
measurements and estimates the sensor positions. Distributed
approaches require that each node computes its own position, working
only with data collected locally and exchanged with neighbors. The
body of work on centralized approaches to the problem is vast (\eg,
\cite{OguzGomesXavierOliveira2011,DestinoAbreu2011,BiswasLiangTohYeWang2006,KhanKarMoura2010,BlattHero2006}). At
a smaller scale, distributed techniques based on convex relaxations of
the problem are also abundant (for example,
\cite{GholamiTetruashviliStromCensor2013,SimonettoLeus2014,KellerGur2011,CostaPatwariHero2006,SrirangarajanTewfikLuo2008})
But distributed, and maximum likelihood (thus nonconvex) approaches to
the sensor network localization problem are much less common. The
algorithm presented in~\cite{ShiHeChenJiang2010} is a nonlinear
Gauss-Seidel approach: only one node works at a time and solves a
source localization problem with neighbors playing the role of
anchors. The nodes activate sequentially in a round-robin
scheme. Thus, the time to complete just one cycle becomes proportional
to the network size. Parallel algorithms --- the ones we are
interested in this paper --- avoid altogether this issue, as all nodes
operate simultaneously; moreover, adding or deleting a node raises no
special synchronization concern. The work presented
in~\cite{CalafioreCarloneWei2010} puts forward a two-stage algorithm
which is parallel: in a first consensus phase, a Barzilai-Borwein (BB)
step size is calculated, followed by a local gradient computation
phase. It is known that BB steps do not necessarily decrease the
objective function; as discussed in \cite{Raydan1997}, an outer
globalization scheme involving line searches is needed to ensure its
stability. However, line searches are cumbersome to implement in a
distributed setting and are, in fact, absent
in~\cite{CalafioreCarloneWei2010}. Further, the algorithm requires the
step size to be computed via consensus, and thus the number of
consensus rounds needed is a parameter to tune. We will present an
algorithm with simple implementation which is both parallel and
stable, with no free parameters.  We will compare experimentally the
performance of our method with the distributed, parallel, state of the
art method in~\cite{CalafioreCarloneWei2010} in
Sec.~\ref{sec:experimental-results}.

\paragraph*{Contributions}
\label{sec:contributions}

We tackle directly the nonconvex problem
in~\eqref{eq:snlOptimizationProblem}, with a simple and efficient
 algorithm which:
\begin{enumerate}
\item is parallel;
\item does not involve any free parameter definition; 
\item is proven not to increase the value of the cost function at each iteration;
\item has better performance in positioning error and cost value than
  the state of the art method, while expending less in
  communications.
\end{enumerate}
The first and second claims are addressed in
Sec.~\ref{sec:distr-sens-netw}, the third in
Sec.~\ref{sec:major-minim} and the last one in
Sec.~\ref{sec:experimental-results}, dedicated to numerical
experiments.

\section{Problem reformulation}
\label{sec:problem-reformulation}

We can reformulate Problem~(\ref{eq:snlOptimizationProblem}) as
\begin{IEEEeqnarray}{ll}
  \nonumber
    \minimize_{x_{i}, y_{ij}, w_{ik}} & \sum_{i \sim j} \frac{1}{2}
    \|x_{i} - x_j - y_{ij}\|^2 + \\  \label{eq:reformulation}
    & \sum_i \sum_{j \in {\cal A}_i} \frac{1}{2} \|x_{i} - a_k -
    w_{ik}\|^2\\ \nonumber
    \text{subject to } & \|y_{ij}\| = d_{ij}, \;\|w_{ij}\| = r_{ij},
\end{IEEEeqnarray}
and rewrite~(\ref{eq:reformulation}) as
\begin{IEEEeqnarray}{ll}
  \label{eq:reformulation2}
    \minimize_{x_{i}, y_{ij}, w_{ik}} & \frac{1}{2} \|Ax - y\|^2 +
    \sum_i \frac{1}{2} \|x_{i} \otimes 1 - \alpha_i - w_i\|^2\\
    \nonumber
    \text{subject to } & \|y_{ij}\| = d_{ij}, \; \|w_{ik}\| = r_{ik},
\end{IEEEeqnarray}
with concatenated vectors $x = (x_{i})_{i \in \mathcal{V}}$, $y =
(y_{ij})_{i \sim j}$, $\alpha_{i}=(a_{ik})_{k \in \mathcal{A}_{i}}$,
and $w_{i} = (w_{ik})_{k \in
  \mathcal{A}_{i}}$. In~\eqref{eq:reformulation2}, the symbol~$1$
stands for the vector of ones. Matrix~$A$ is the
result of the Kronecker product of the arc-node incidence
matrix\footnote{Each edge is arbitrarily assigned a direction by the two
  incident nodes.}~$C$
with the identity matrix~$I_{p}$: $A = C \otimes I_{p}$.
Problem~\eqref{eq:reformulation2} is equivalent to
\begin{IEEEeqnarray*}{ll}
  \minimize_{x_{i}, y_{ij}, w_{ik}} & \frac{1}{2} \left \|
    \begin{bmatrix}
      A & -I & 0
    \end{bmatrix}
    \begin{bmatrix}
      x\\
      y\\
      w
    \end{bmatrix}
  \right \|^{2} + \frac{1}{2} \left \| Ex - \alpha -
    w\right\|^{2}\\
  \text{subject to } & \|y_{ij}\| = d_{ij}, \;\|w_{ik}\| = r_{ik},
\end{IEEEeqnarray*}
where $\alpha = (\alpha_{i})_{i \in \mathcal{V}}$, $w = (w_{i})_{i \in
  \mathcal{V}}$, and $E$ is a matrix with zeros and ones, selecting
the entries in $\alpha$ and $w$ corresponding to each sensor node. We
now collect all the optimization variables in~$z =
(x,y,w)$, and rewrite our problem as
\begin{IEEEeqnarray*}{ll}
  \label{eq:reform3}
    \minimize_{z} & \frac{1}{2} \left \|
      \begin{bmatrix}
        A & -I & 0
      \end{bmatrix}
      z \right \|^{2} + \frac{1}{2} \left \|
      \begin{bmatrix}
        E & 0 & -I
      \end{bmatrix}
      z - \alpha \right \|^{2}\\
    \text{subject to } & z \in \mathcal{Z},
\end{IEEEeqnarray*}
where $\mathcal{Z} = \{ z = (x, y, w): \|y_{ij}\| = d_{ij}, i \sim j,
w_{ik} = r_{ik}, i \in \mathcal{V}, k \in
\mathcal{A}_{i}\}$. Problem~\eqref{eq:reformulation2} can be written as
\begin{IEEEeqnarray}{ll}
  \label{eq:reformulation4}
      \minimize_{z} & f(z) = \frac{1}{2} z^{T} M z - b^{T} z \\
    \text{subject to } & z \in \mathcal{Z},
  \end{IEEEeqnarray}
for $M$ and $b$ defined as
\begin{IEEEeqnarray}{rClrCl}
\label{eq:M1-M2}
  M &=& M_{1} + M_{2},& \quad b &=&
  \begin{bmatrix}
    E^{T}\\0\\-I
  \end{bmatrix}
  \alpha,
\end{IEEEeqnarray}
\begin{IEEEeqnarray*}{rClrCl}
  M_{1} &=& 
\begin{bmatrix}
    A^{T}\\-I\\0
  \end{bmatrix}
  \begin{bmatrix}
    A & -I & 0
  \end{bmatrix},
  & \quad M_{2} &=& 
  \begin{bmatrix}
    E^{T}\\ 0\\ -I
  \end{bmatrix}
  \begin{bmatrix}
    E & 0 & -I
  \end{bmatrix}.
\end{IEEEeqnarray*}

\section{Majorization-Minimization}
\label{sec:major-minim}

To solve Problem~\eqref{eq:reformulation4} in a distributed way we
must deal with the complicating off-diagonal entries of~$M$ that
couple the sensors' variables. We emphasize a simple, but key fact:
\begin{remark}
  \label{th:remark-quadratic-lipschitz}
  The function optimized in Problem~\eqref{eq:reformulation4} is
  quadratic in $z$ and, thus, has a Lipschitz continuous
  gradient~\cite{bertsekas1999nonlinear},
  \ie,
  \begin{equation*}
    \|\nabla f(x) -\nabla f(y)\| \leq L\|x-y\|,
  \end{equation*}
  for some~$L$ and all $x,y$.
\end{remark}
From this property of function~$f$ we can obtain the upper bound (also
found in~\cite{bertsekas1999nonlinear}) $ f(z) \leq f(z^{t}) + \left <
  \nabla f(z^{t}), z - z^{t} \right> + \frac{L}{2} \left \|z - z^{t}
\right\|^{2}$, for any point~$z^{t}$ and use it as a majorizer in the
Majorization-Minimization framework\cite{HunterLange2004}. This
majorizer decouples the variables and allows for a distributed
solution. Our algorithm is simply:
\begin{equation}
  \label{eq:iterzArgmin}
  z^{t+1} = \argmin_{z \in \mathcal{Z}} f(z^{t}) + \left < \nabla
    f(z^{t}), z - z^{t} \right> +
  \frac{L}{2} \left \|z - z^{t} \right\|^{2}.
\end{equation}
The solution of~\eqref{eq:iterzArgmin} is the projected gradient
iteration~\cite{bertsekas1999nonlinear}
\begin{equation}
  \label{eq:proj-grad}
  z^{t+1} = \mathrm{P}_{\mathcal{Z}}\left(z^{t} - \frac{1}{L}\nabla
    f(z^{t}) \right),
\end{equation}
where~$\mathrm{P}_{\mathcal{Z}}(p)$ is the projection of point~$p$
onto~$\mathcal{Z}$. The gradient in~\eqref{eq:proj-grad} can be easily
computed as the affine function
$\nabla f(z) = Mz-b.$
See the recent
work~\cite{BeckEldar2013} for interesting convergence properties of
the recursion~~\eqref{eq:proj-grad}. Particularly, we emphasize that
the cost function is non increasing per iteration.

We now compute a Lipschitz constant~$L$ for the gradient of the
quadratic function in Problem~\eqref{eq:reformulation4}, such that it
is easy to estimate in a distributed way.
\begin{IEEEeqnarray}{rCl}
  \nonumber L & = & \lambda_{\mathrm{max}}\left(M\right)\\ \nonumber
  &\leq& \lambda_{\mathrm{max}}\left(M_{1}\right) +
  \lambda_{\mathrm{max}}\left(M_{2}\right) \\ \nonumber
  & = & \lambda_{\mathrm{max}}\left(AA^{T}+I\right) +
  \lambda_{\mathrm{max}}\left( EE^{T} + I\right) \\ \nonumber
  &\leq& \lambda_{\mathrm{max}}\left(A^{T}A \right) +
  \lambda_{\mathrm{max}}\left( EE^{T}\right) + 2 \\
  \label{eq:lips-const}
  &\leq&2\delta_{\mathrm{max}} + \max_{i \in \mathcal{V}}
  |\mathcal{A}_{i}| + 2,
\end{IEEEeqnarray}
where~$\lambda_{\mathrm{max}}$ denotes the largest eigenvalue,
$|\mathcal{A}|$ is the cardinality of set~$\mathcal{A}$,
and~$\delta_{\mathrm{max}}$ is the maximum node degree of the
network. We note that $\lambda_{\mathrm{max}}(A^{T}A)$ is the maximum
eigenvalue of graph~$\mathcal{G}$ laplacian matrix; the proof that it
is upper-bounded by~$2\delta_{\mathrm{max}}$ can be found in~\cite{chung1997spectral}. This Lipschitz constant can be computed
in a distributed way by, \eg, a diffusion algorithm
(c.f. \cite[Ch.~9]{MesbahiEgerstedt2010}).

\section{Distributed sensor network localization}
\label{sec:distr-sens-netw}

At this point, the recursion in Eq.~\eqref{eq:proj-grad} is
already distributed, as detailed below.
From~\eqref{eq:proj-grad} we will obtain the update rules for
the variables~$x$,~$y$ and~$w$. For this we write matrix~$M$ as follows:
\begin{equation}
  \label{eq:M-and-b-per-var}
  M = 
  \begin{bmatrix}
    A^{T}A + E^{T}E & -A^{T} & -E^{T}\\
    -A & I & 0\\
    -E & 0 & I
  \end{bmatrix},
\end{equation}
and denote~$B = A^{T}A + E^{T}E$. Then, each block of~$z$ is
updated according to
\begin{algorithm}[tb]
  \caption{Distributed nonconvex localization algorithm}
  \label{alg:distr-loc}
  \begin{algorithmic}[1]
    \REQUIRE $x^{0}; L; \{d_{ij} : j \in N_{i}\}; \{r_{ik} : k
    \in \mathcal{A}_{i}\};$
    \ENSURE $\hat{x}$
    \STATE set $y_{ij}^{0} =
    \mathrm{P}_{\mathcal{Y}_{ij}}\left(x^{0}_{i} -x^{0}_{j}\right)$,
    $
    \mathcal{Y}_{ij} = \{y: \|y\| = d_{ij}\}
    $
    and $w_{ik}^{0} =
    \mathrm{P}_{\mathcal{W}_{ik}}\left(x^{0}_{i} -a_{k}\right)$,
    $
    \mathcal{W}_{ik} = \{w: \|w\| = r_{ik}\}
    $
    \STATE $t = 0$
    \WHILE{some stopping criterion is not met, each node~$i$}
    \STATE $
          x_{i}^{t+1} =
      b_{i}x_{i}^{t} +
      \frac{1}{L}\sum_{j \in N_{i}}\left(x_{j}^{t} + C_{(i \sim
          j,i)}y_{ij}^{t}\right) + 
      \frac{1}{L} \sum_{k \in
        \mathcal{A}_{i}} \left(w_{ik}^{t} + a_{ik}\right)
     $\label{alg:x}
    \STATE for all neighboring~$j$, compute\\
    $
      y_{ij}^{k+1} = \mathrm{P}_{\mathcal{Y}_{ij}} \left(
        \frac{L-1}{L}y_{ij}^{k} + \frac{1}{L} C_{(i \sim
          j,i)}\left(x_{i}^{t}-x_{j}^{t}\right)\right),
    $
    \STATE for each of the connected anchors~$k \in \mathcal{A}_{i}$, compute\\
    $
      w_{ik}^{t+1} = \mathrm{P}_{\mathcal{W}_{ik}}\left(
        \frac{L-1}{L}w_{ik}^{t}+\frac{1}{L}(x_{i}-a_{ik})\right)
    $
    \STATE broadcast~$x_{i}^{t+1}$ to neighbors
    \STATE $t = t + 1$
    \ENDWHILE
    \RETURN $\hat x_{i} = x_{i}^{t}$
  \end{algorithmic}
\end{algorithm}
\begin{IEEEeqnarray}{rCl}
  \label{eq:x-update-rule}
  x^{t+1} &=& \left( I - \frac{1}{L}B\right) x^{t} +
  \frac{1}{L}A^{T}y^{t} + \frac{1}{L}E^{T}(w^{t} + \alpha),\\
  \label{eq:y-update-rule}
  y^{t+1} &=& \mathrm{P}_{\mathcal{Y}}\left( \frac{L-1}{L}y^{t} +
    \frac{1}{L}Ax^{t}\right),\\
  \label{eq:w-update-rule}
  w^{t+1} &=& \mathrm{P}_{\mathcal{W}}\left( \frac{L-1}{L} w^{t} +
    \frac{1}{L}Ex^{t} - \frac{\alpha}{L}\right),
\end{IEEEeqnarray}
where~$\mathcal{Y}$ and~$\mathcal{W}$ are the constraint sets
associated with the acquired measurements between sensors, and between
anchors and sensors, respectively, and~$N_{i}$ is the set of the
neighbors of node~$i$. We observe that each block of~$z = (x, y, w)$
at iteration~$t+1$ will only need local neighborhood information, as
clarified in Algorithm~\ref{alg:distr-loc}. Each node~$i$ will update
the current estimate of its own position, each one of the~$y_{ij}$ for
all the incident edges~$i \sim j$ and the anchor terms~$w_{ik}$, if
any. The symbol $C_{(i \sim j,i)}$ denotes the arc-node incidence
matrix entry relative to edge~$i \sim j$ (row index) and node~$i$
(column index). The constant~$b_{i}$ in step~\ref{alg:x} of
Algorithm~\ref{alg:distr-loc} is defined
as~$\frac{L-\delta_{i}-|\mathcal{A}_{i}|}{L}$.

\section{Experimental results}
\label{sec:experimental-results}

We present numerical experiments to ascertain the performance of the proposed
Algorithm~\ref{alg:distr-loc}, both in accuracy and in communication
cost. Accuracy will be measured in 1)~\emph{mean positioning error per sensor},
defined as
\begin{equation}
  \label{eq:mean-error}
  MPE = \frac{1}{n \cdot MC} \sum_{mc=1}^{MC}\sum_{i=1}^{n}
  \|\hat{x}_{i}(mc) - x_{i}^{\star}\|,
\end{equation}
where~$MC$ is the total number of Monte Carlo
trials,~$\hat{x}_{i}(mc)$ is the estimate generated by an algorithm at
the Monte Carlo trial~$mc$, and~$x_{i}^{\star}$ is the true position
of node~$i$, and 2)~also by evaluating the cost function
in~\eqref{eq:snlOptimizationProblem}, averaged by the Monte Carlo
trials and number of sensors, as
in~\eqref{eq:mean-error}. Communication cost will be measured taking
into account that each iteration in Algorithm~\ref{alg:distr-loc}
involves communicating~$pn$ real numbers. We will compare the
performance of the proposed method with the Barzilai-Borwein algorithm
in~\cite{CalafioreCarloneWei2010}, whose communication cost per
iteration is $n(2T + p)$, where~$T$ is the number of consensus rounds
needed to estimate the Barzilai-Borwein step size. We use $T = 20$ as
in~\cite{CalafioreCarloneWei2010}. The setup for the experiments is a
geometric network with $50$~sensors randomly distributed in the
two-dimensional square~$[0, \: 1] \times [0, \: 1]$, with average node
degree of about~$6$, and $4$~anchors placed at the vertexes of this
square. The network remains fixed during all the Monte Carlo trials.
Both algorithms receive an initialization from a convex approximation
method. The initialization will hopefully hand the nonconvex
refinement algorithms a point near the basin of attraction of the true
minimum. For this purpose we generate noisy range measurements
according to $d_{ij} = | \|x_{i}^{\star} - x_{j}^{\star}\| + \nu_{ij}
|$, and $r_{ik} = | \|x_{i}^{\star} - a_{k}\| + \eta_{ik} |$, where
~$\{\nu_{ij} : i \sim j \in \mathcal{E}\} \cup \{\eta_{ik} : i \in
\mathcal{V}, k \in \mathcal{A}_{i}\}$ are independent gaussian random
variables with zero mean and standard deviation~$\sigma$.  We
conducted $100$~Monte Carlo trials for each standard
deviation~$\sigma=(0.01, \: 0.05, \: 0.1)$. If we spread the sensors
by a squared area with side of $1$Km, this means measurements are
affected by noise of standard deviation of $10$m, $50$m, and $100$m.
\begin{table}[tb]
  \centering
  \caption{Mean positioning error, with measurement noise}
  \label{tab:errVsStdev}
  \begin{tabular}{@{}rrr@{}}
    \toprule
    \textbf{$\sigma$}&\textbf{Proposed
      method}&\textbf{BB method}\\\midrule
    0.01&0.0053&0.0059\\
    0.05&0.0143&0.0154\\
    0.10&0.0210&0.0221\\
    \bottomrule
  \end{tabular}
\end{table}
In terms of mean positioning error per sensor the proposed algorithm
fares better than the benchmark: Table~\ref{tab:errVsStdev} shows the
mean error defined in~\eqref{eq:mean-error} after the algorithms have
stabilized, or reached a maximum iteration number. In a square
with~$1$Km sides, we improve the accuracy of the gradient descent with
Barzilai-Borwein steps by about $1$m per sensor, even for high power
noise.
\begin{figure}[!h]
  \centering \subfigure[The proposed method improves the comparing
  algorithm, both in accuracy and communication cost. Our proposed method
  improves the state of the art method
  in~\cite{CalafioreCarloneWei2010} by about $60$~cm in mean positioning
  error per sensor, delivering a no surprises, stable progression of
  the error of the estimates.]
  { \label{fig:errVsCommunications-0_01}
\includegraphics[width=0.95\columnwidth]{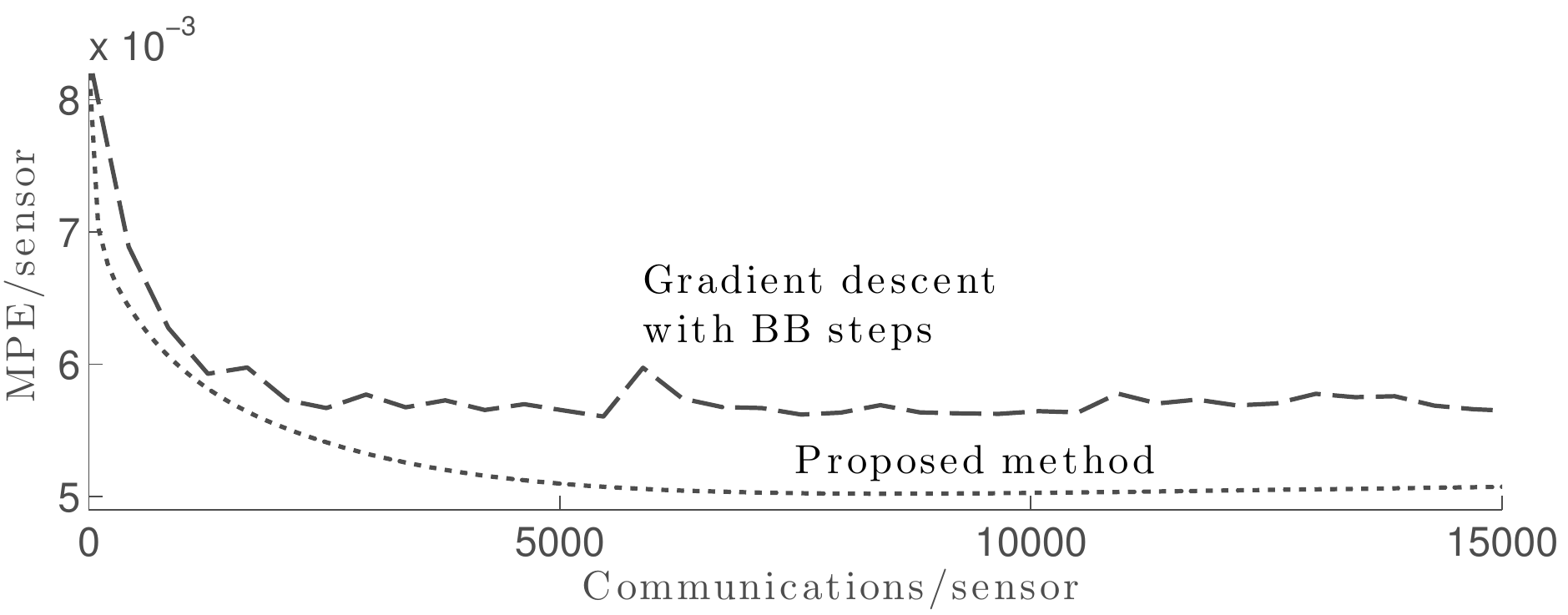}
   }
  \subfigure[The final
  costs are, for the BB method,~$1.7392\: 10^{-4}$ and, for the proposed
  method~$1.5698\: 10^{-4}$. A small difference in cost that
  translates into a considerable distance in error, as depicted in
  Fig.~\ref{fig:errVsCommunications-0_01} and
  Table~\ref{tab:errVsStdev}.]
  {\label{fig:costVsCommunications-0_01}
    \includegraphics[width=0.95\columnwidth]{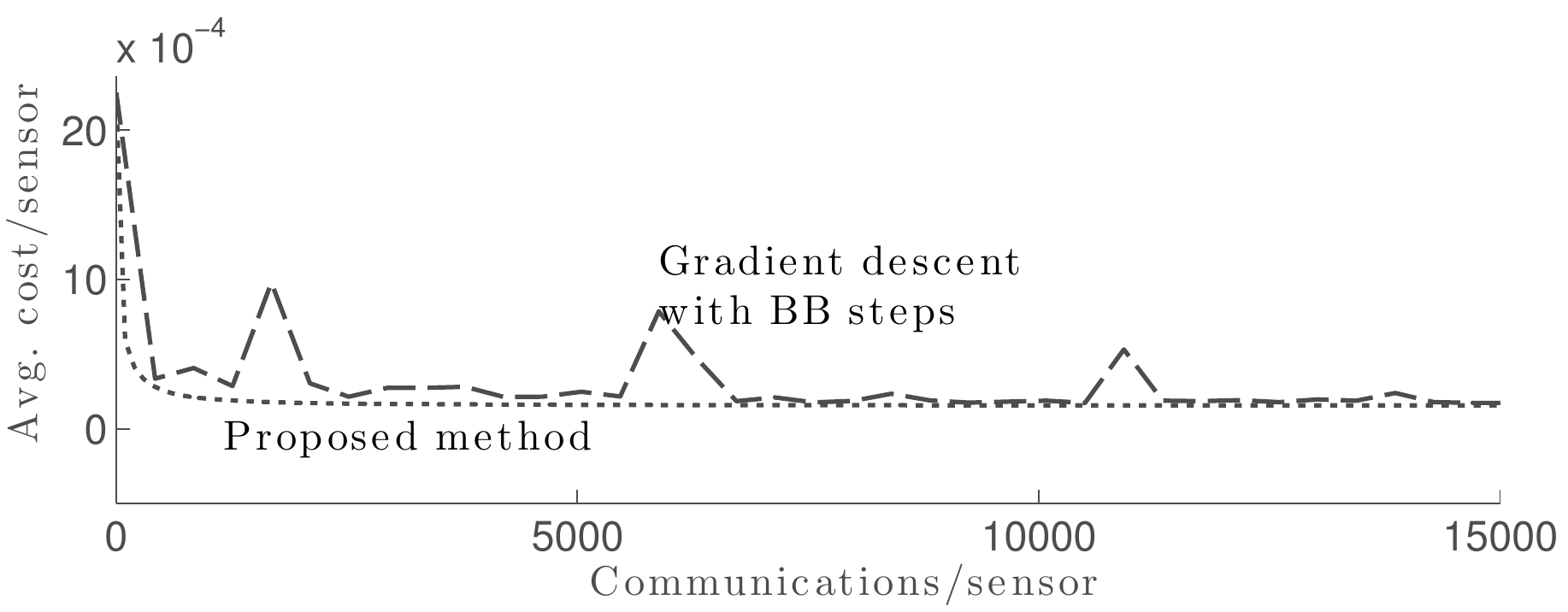}
  } 
  \caption{Noisy distance measurements with~$\sigma = 0.01$,
    representing~$10$m for a square with $1$Km sides. The proposed
    method shows a faster and smoother progression, while the
    comparing method bounces, always above the proposed method.}
  \label{fig:noise-0_01}
\end{figure}
Fig.~\ref{fig:noise-0_01} depicts the averaged evolution of the
error per sensor of both algorithms as a function of the volume of
accumulated communications, and also the evolution of the cost. The
gradient descent with Barzilai-Borwein steps shows an irregular
pattern for the error, only vaguely matching the variation in the
corresponding cost (Fig.~\ref{fig:costVsCommunications-0_01}), thus
leaving some uncertainty on when to stop the algorithm and what
estimate to keep.
The presented
method reaches the final cost value per sensor much faster and
steadily than the benchmark for medium-low measurement
noise. In fact, our method takes under one order of magnitude less
communications than the comparing one to approach the minimum cost
value (match the cost at about~$1500$ communications with~$15000$).
\begin{figure}[!h]
  \centering 
  \subfigure[For medium noise power the algorithms' performance
  comparison follows the one under low noise power. The accuracy gain
  is more than~$1$m per sensor.]
  {\label{fig:errVsCommunications-0_05}
    \includegraphics[width=0.95\columnwidth]{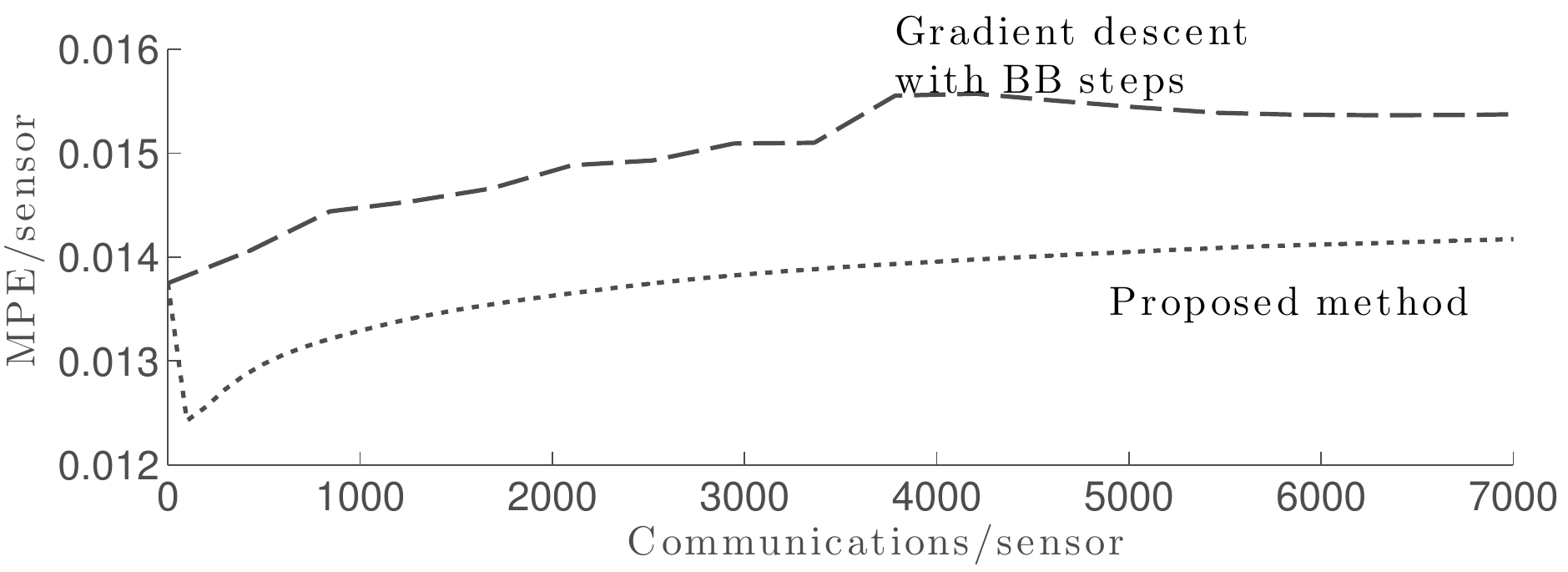}}
  \subfigure[Under medium noise the proposed method also reaches a
  smaller value for the average cost per sensor:~$0.0031$,
  against~$0.0032$ from the BB method.]
  {\label{fig:costVsCommunications-0_05}
    \includegraphics[width=0.95\columnwidth]{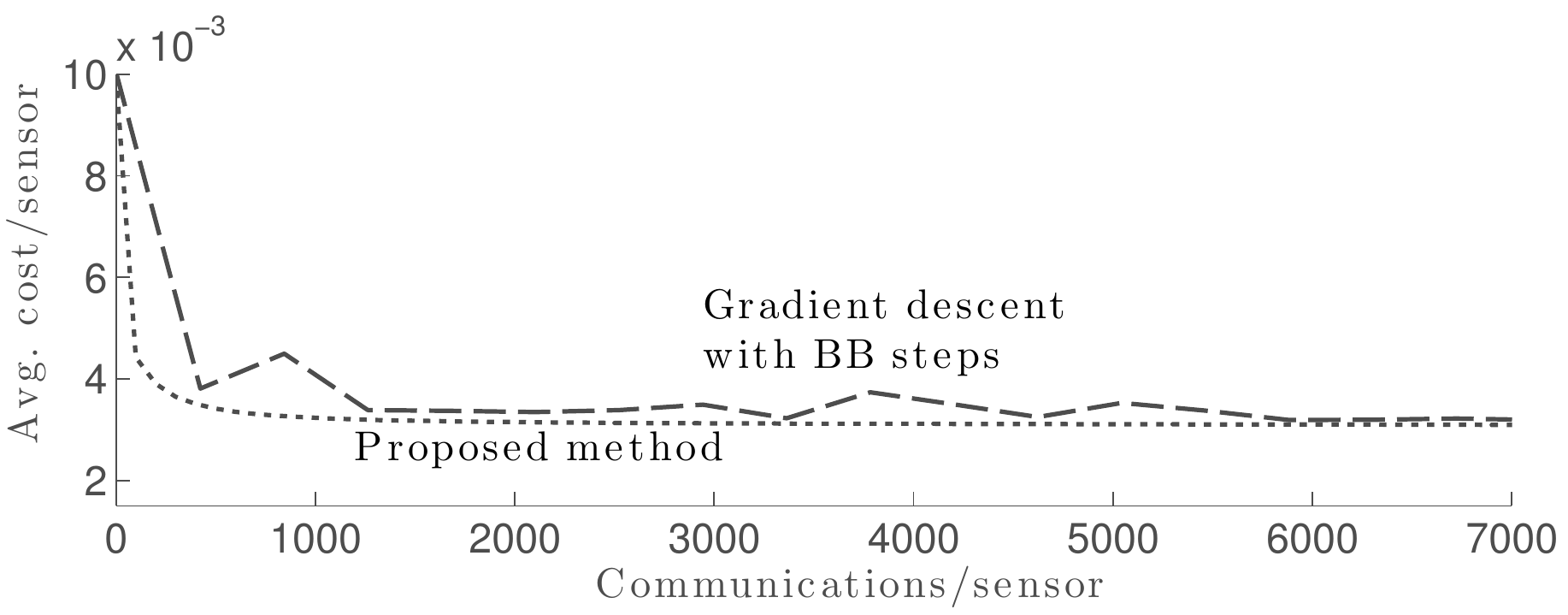}
  }
  \caption{Distance measurements contaminated with noise,
    with~$\sigma = 0.05$, representing~$50$m for a square with $1$Km
    sides. The proposed method continues to outperform the
    comparing state of the art method, and contrasting the
    instability of the BB method.}
\label{fig:noise-0_05}
\end{figure}
The most realistic case of medium noise power led to the results
presented in Fig.~\ref{fig:noise-0_05}. The characteristic
irregularity of the BB method continues to fail in delivering better
solutions in average than our stable, guaranteed method. The error
curves in Fig.~\ref{fig:errVsCommunications-0_05} are increasing,
because the error is not the quantity being directly optimized and the
medium-high noise power in measurement data shifts the cost optimal
points.
\begin{figure}[!h]
  \centering \subfigure[The proposed algorithm tops the comparing
  method in error, under high noise power, by more than~$1$m, when
  considering a squared deployment area of~$1$Km sides.]
{\includegraphics[width=0.95\columnwidth]{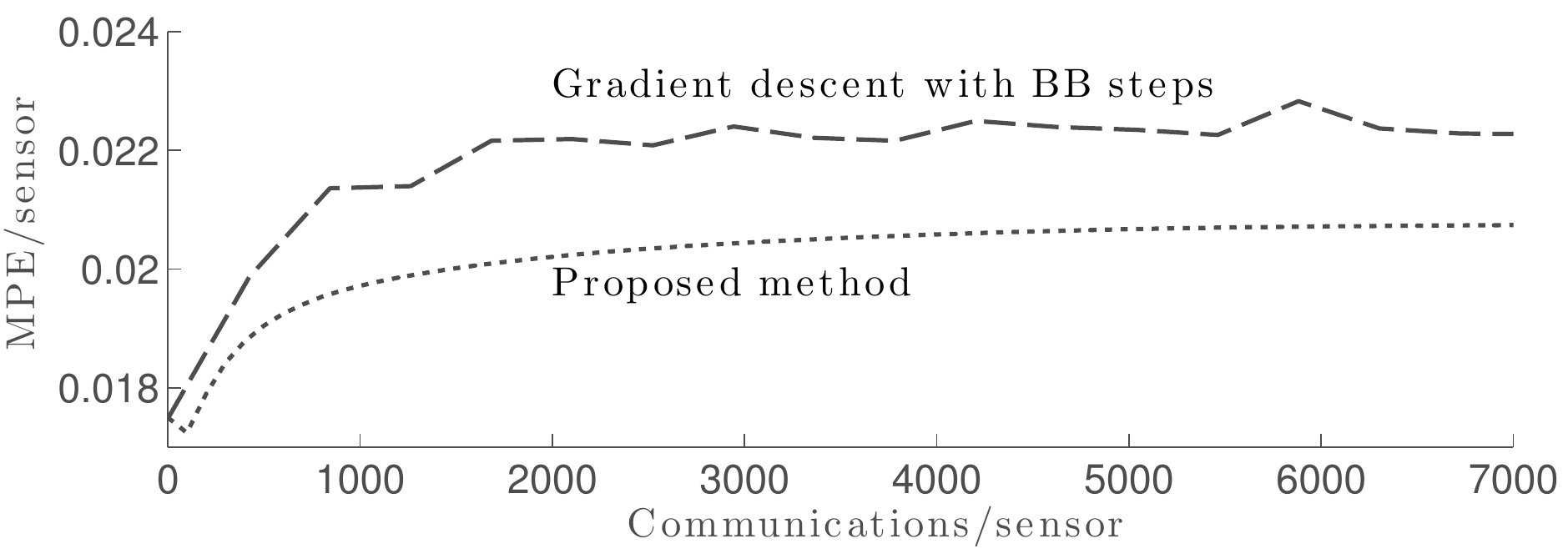}
    \label{fig:errVsCommunications-0_1}}
  \subfigure[Under heavy noise the proposed method reaches a smaller value
  for the average cost per sensor:~$0.0096$, against~$0.0099$ from
  the BB method.]
  {\includegraphics[width=0.95\columnwidth]{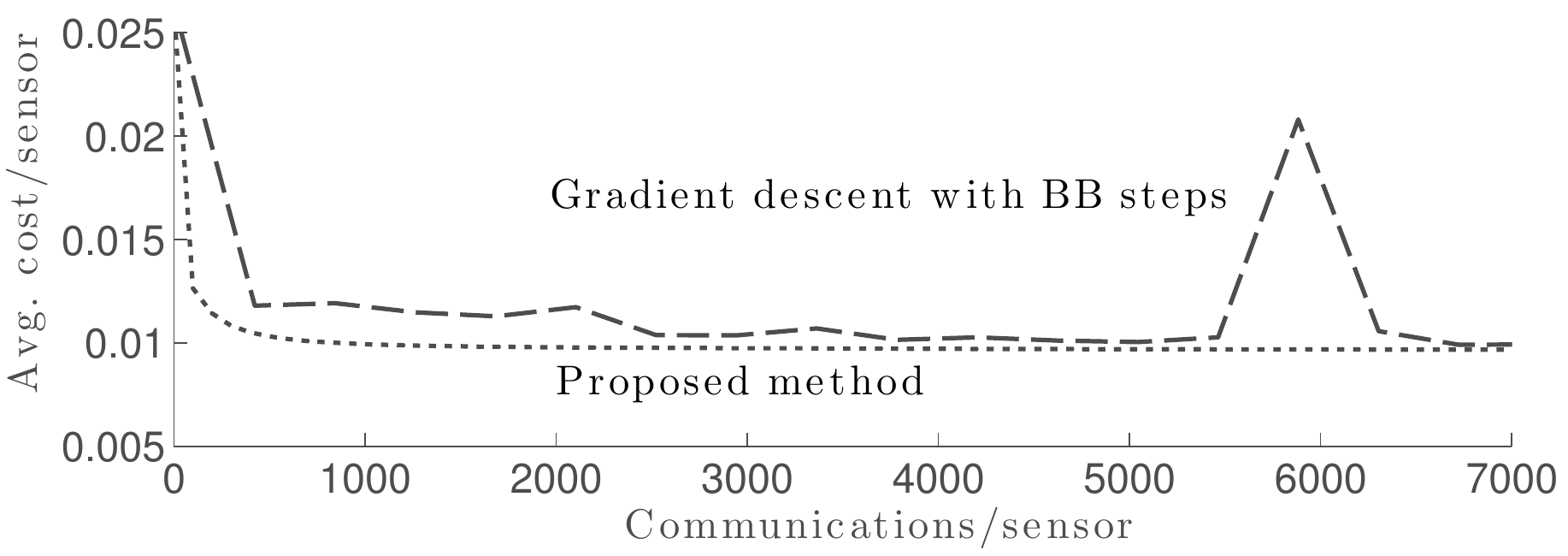}
    \label{fig:costVsCommunications-0_1}}
\label{fig:noise-0_1}
\caption{Distance measurements contaminated with noise, with~$\sigma =
  0.1$, representing~$100$m for a square with $1$Km sides.}
\end{figure}
Under high noise power, our method tops the performance of the
benchmark in cost function terms, as it is shown in
Fig.~\ref{fig:costVsCommunications-0_1}, not only in the convergence
speed, but also in the final value reached. Again, our method has
almost one order of magnitude less in communications to achieve its
plateau, which is itself, in average, better than the alternative method
(compare the performance at~$700$ communications with the one at~$7000$).

\section{Concluding remarks}
\label{sec:concluding-remarks}

The monotonicity of the proposed method is a strong feature for 
applications of sensor network localization.  Our method proves to be
not only fast and resilient, but also simple to implement and deploy,
with no free parameters to tune. The steady accuracy gain over the
competing method also makes it usable in contexts with a wide range of
measurement errors are expected. The presented method can be useful
both as a refinement algorithm and as a tracking method, \eg, for mobile
robot formations where position estimates computed on a given time
step are used as initialization for the next one.

\bibliographystyle{IEEEtran} \bibliography{IEEEabrv,biblos.bib}
\end{document}